\documentclass{amsart}
\usepackage{graphicx}
\usepackage{graphicx}
\usepackage{amssymb,amscd,amsthm,amsxtra}
\usepackage{latexsym}
\usepackage{epsfig}

\newtheorem{thm}{Theorem}[section]

\newtheorem{lem}[thm]{Lemma}

\theoremstyle{definition}

\theoremstyle{remark}

\numberwithin{equation}{section}
\newcommand{\R}{\mathbb R}

\begin{document}

\title{Bernstein-type techniques for 2D free boundary graphs}%
\author{Daniela De Silva}
\address{Department of Mathematics, Johns Hopkins University, Baltimore, MD 21218}%
\email{desilva@math.jhu.edu}

%\thanks{}%
%\subjclass{}%
%\keywords{}%

%\date{}%
%\dedicatory{}%
%\commby{}%
% ----------------------------------------------------------------
\begin{abstract} We prove an a-priori bound for the Lipschitz
constant of a smooth one-phase free boundary graph $F(u)$ in two
dimensions. The function $u$ satisfies an elliptic equation in its
positive side, and $|\nabla u|=1$ on $F(u).$
\end{abstract}
\maketitle
% ----------------------------------------------------------------
\section{Introduction}

We consider the following one-phase free boundary problem,

\begin{equation}\label{fbintro} \left \{
\begin{array}{ll}
    \Delta u =0,   & \hbox{in $\Omega^+(u):= \{x \in \Omega : u(x)>0\}$,} \\
    |\nabla u| =1, & \hbox{on $F(u):= \partial \Omega^+(u) \cap \Omega$,} \\
\end{array}\right.
\end{equation}
\vspace{0.01in}

\noindent where $\Omega$ is a domain in $\R^n.$ Here, for any
non-negative function $v : \Omega \rightarrow \mathbb{R}$, we set

$$\Omega^+(v) := \{x \in \Omega : v(x) >0\}, \ \ \Omega^-(v):=\{x \in \Omega :
v(x)=0\}^\circ, \ \ F(v) := \partial \Omega^+(v) \cap \Omega.$$
\vspace{0.01in}

\noindent Problem \eqref{fbintro} arises for example in the
minimization of the variational integral

$$J(u)=
\int_{\Omega}(|\nabla u|^2 + \chi_{\{u>0\}})dx,$$ \vspace{0.01in}

\noindent that appears in many applications (see \cite{AC},
\cite{F}.)

In \cite{C1},\cite{C2},\cite{C3}, the author introduced the notion
of ``viscosity" solution to \eqref{fbintro}, and developed the
theory of existence and regularity of viscosity free boundaries.
In particular, the regularity theory is inspired by the regularity
theory for minimal surfaces, precisely by the ``oscillation decay"
method of De Giorgi, according to which if $S$ is a minimal
surface in the unit ball $B_1$, and $S$ is the graph of a
Lipschitz function, then $S$ is $C^{1,\alpha}$ (hence smooth) in
$B_{1/2}$. Analogously, if $F(u)$ is a Lipschitz free boundary in
$B_1$, then $F(u)$ is $C^{1,\alpha}$ in $B_{1/2}.$ Higher
regularity results of \cite{KN} then yield the local analyticity
of $F(u)$ in the interior.

Thus, a natural question arises, that is how to obtain the
Lipschitz continuity of a viscosity free boundary. In the theory
of minimal surfaces, in the special case of minimal graphs, this
is achieved via an a-priori gradient bound for solutions to the
minimal surface equation, originally proved in \cite{BMG}.
Analogously, an a-priori bound for the Lipschitz constant of
smooth free boundary graphs is needed, in order to obtain that
viscosity free boundary graphs are smooth in the interior.

In this note we provide this tool in the 2D and 3D case. Our proof
is based on the so-called Bernstein technique, which is been
widely used in literature (see for example \cite{GT}.) A similar
approach for minimal surfaces is used in \cite{WX}.

Moreover, in the 2D case, our technique is flexible enough to
allow us to obtain an a-priori bound for the Lipschitz constant of
free boundary graphs for a wider class of problems, to which the
regularity theory of \cite{C1} has been extended (see for example
\cite{WP}).

In order to state our main result, we introduce some notation. Let

$$C_{(R,M)} := \mathcal{B}_R(0) \times [-M,M] \subset
\mathbb{R}^n.$$\vspace{0.01in}

\noindent Here $\mathcal{B}_r$ denotes a $(n-1)$-dimensional ball,
while a $n$-dimensional ball is denoted by $B_r$. When $R=1$, we
simply write $C_M$ for $C_{(1,M)}.$ Also, a point $x \in
\mathbb{R}^n$ is denoted by $(x',x_n).$

Assume that $u$ is a classical solution to the following one-phase
free boundary problem:

\begin{equation}\label{fbproblem} \left \{
\begin{array}{ll}
     \mathcal{F}(D^2u)=0,   & \hbox{in $C^+_{M}(u)$,} \\
    |\nabla u| =1, & \hbox{on $F(u)$,} \\
   u_n >0 , & \hbox{on $\overline{C^+_{M}(u)}$,} \\
\end{array}\right.
\end{equation}\vspace{0.01in}

\noindent that is, $u \in C^2(\overline{C^+_{M}(u)})$, and $F(u)$
is a $C^2$-surface. Here $\mathcal{F}$ is a nonlinear uniformly
elliptic operator with ellipticity constants $0< \lambda \leq
\Lambda,$ and $\mathcal{F}(0)=0.$

\noindent Furthermore, assume that $0 \in F(u)$, and that

$$B_{2\alpha }(0, M/2) \Subset
 C^+_{M}(u), \ \ B_{2\alpha }(0, -M/2) \Subset C^-_{M}(u),$$ \vspace{0.01in}

\noindent for some constant $\alpha <1/2.$

By the implicit function theorem, $F(u)$ is a smooth graph in the
$x_n$ (vertical) direction. Let us denote by $Lip(s)$ the
Lipschitz constant of $F(u)$ over $B_s(0)$.

In this note we focus on the 2 dimensional case, $n=2.$ Our main
result is the following a priori bound.

\begin{thm}\label{main} Assume $n=2,$ and let $u$ be a solution to \eqref{fbproblem}. Then, there exists constants C, s,
depending on $M, \lambda, \Lambda, \alpha$ such that

$$\frac{|u_1|}{u_2}\leq C \ \ \textit{on} \ \ B_s(0).$$\vspace{0.01in}

\noindent In particular $Lip(s) \leq C.$
\end{thm}

The paper is organized as follows. In Section 2, we prove a
technical Lemma, which is dimension independent. In Section 3, for
expository purposes, we present the proof of Theorem \ref{main} in
the case when $\mathcal{F}(D^2u) = \Delta u.$ Then, in Section 4,
we proceed to show the proof of Theorem \ref{main} in the general
case. In Section 5, we obtain the analogue of Theorem \ref{main}
in the case when $u$ is a $p$-harmonic function in its positive
phase. Finally, we conclude this note with some remarks about the
applicability of the method of the proof to problems in higher
dimensions.

\section{A preliminary Lemma.}

Here and henceforth, $C, C'$ will denote constants depending
possibly on $M, \lambda, \Lambda$ and $\alpha.$

We start with the following technical Lemma, which holds in any
dimension.

\begin{lem}\label{bound} Let $u$ be a solution to \eqref{fbproblem}. Then, there exist constants $C,C'$, such that
\begin{enumerate}\item$|\nabla u| \leq C'$ on
$C_{(\alpha,M/2)},$\item $u \geq C$ on
$B_{\alpha}(0,M/2)$.\end{enumerate}
\end{lem}
\begin{proof}

(1) We start by showing that if $x_0 \in C_{(\alpha,M/2)}^+(u),$
and $d=\text{dist}(x_0,F(u))$, then

\begin{equation}\label{usmall}u(x_0) \leq Cd.\end{equation}Let

$$v(x)=\frac{1}{d}u(x_0 +dx)$$ \vspace{0.01in}

\noindent be the rescale of $u$ in $B_d(x_0)$. Then $v \geq 0$
solves a uniformly elliptic equation

$$\mathcal{G}(D^2v)=0 \ \ \text{in} \ \ B_1(0),$$ \vspace{0.01in}

\noindent with $\mathcal{G}$ having the same ellipticity constants
as $\mathcal{F},$ and $\mathcal{G}(0)=0.$ Hence, by Harnack's
inequality (see \cite{CC})  $$v \geq c v(0) \ \ \text{in} \ \
B_{1/2}(0).$$ \vspace{0.01in}

\noindent Let us choose $\beta<0$ such that, the radially
symmetric function

$$g(x)=\frac{c v(0)}{2^{-\beta} -1}(|x|^\beta -1)$$\vspace{0.01in}

\noindent satisfies $\mathcal{G}(D^2g) \geq 0$ in the anullus $B_1
\setminus \overline{B_{1/2}}$, $g = 0$ on $\partial B_1$ and $g =
c v(0)$ on $\partial B_{1/2}.$ Then, by the maximum principle

$$v \geq g \ \
\text{in} \ \ B_1 \setminus \overline{B_{1/2}}.$$\vspace{0.01in}

\noindent Now, let $x_1 \in \partial B_1(0)$ be such that
$v(x_1)=0$. Then, since $\nabla v(x)= \nabla u(x_0+dx)$, and $u$
solves \eqref{fbproblem}, we have $|\nabla v(x_1)| = 1.$ Let $\nu$
be the inward normal to $\partial B_1$ at $x_1$. Then, at $x_1,$

$$1 = |\nabla v(x_1)| \geq v_{\nu} \geq g_{\nu} \geq C v(0),$$\vspace{0.01in}

\noindent which yields \eqref{usmall}. Using Harnack's inequality
and elliptic regularity (see \cite{CC}) we obtain the desired
claim.

\

(2) Let $g$ be a radially symmetric function such that
$\mathcal{F}(D^2g) \leq 0$, in the annulus $B_{\alpha}(0,-M/2)
\setminus B_{\alpha/2}(0,-M/2)$ and:

$$
  \begin{cases}
    g= C & \text{on $\partial B_{\alpha}(0,-M/2)$}, \\
    g= 0 & \text{on $
 \partial B_{\alpha/2}(0,-M/2)$},\\
    |\nabla g| < 1 & \text{on $\partial B_{\alpha}(0,-M/2).$}
  \end{cases}
$$\vspace{0.01in}

\noindent Since $\mathcal{F} $ is uniformly elliptic and
$\mathcal{F}(0)=0,$ such kind of supersolution can be obtained by
a similar formula as in part $(1).$

Since $u = 0 $ on $B_{2\alpha}(0,-M/2)$, then $g > u$ on the
annulus between $B_{\alpha}(0,-M/2)$ and $B_{\alpha/2}(0,-M/2)$.
Let $g_t$ be a family of translates of $g$ in the positive
vertical direction. Then, the first touching point $x_0$ of $g_t$
and $u$ must occur where $g_t=C$, and $|x'_0| \leq \alpha$, within
distance $\alpha/2$ from $F(u)$. Moreover, $x_0$ occurs before
$F(g_t)$ coincides with $\partial B_{\alpha}(0,M/2)$. Therefore,
since $u$ is monotone increasing in the vertical direction, $u
\geq C$ at some point $x \in \partial B_{\alpha}(0,M/2)$.
Harnack's inequality implies the desired statement, $u \geq C$ on
$B_{\alpha}(0,M/2).$
\end{proof}

\section{Proof of Theorem \ref{main}. The Laplace operator.}

Let us define a smooth positive function $g(x_1,u)$ over the
trapezoid $T(h,a)=\{(x_1,u)\in \mathbb{R}^2 : 0 < u < h, -u-a <
x_1 < u+a \}$, with $h,a>0$, which satisfies the following
properties on $\partial T$:

\begin{enumerate}\renewcommand{\theenumi}{\roman{enumi}}
\item $g(x_1,u) > 0,$ on $B:=\{(x_1,u) : u=0, -a < x_1 < a\};$
\item $g(x_1,u) = 0,$ on $\partial T \setminus B;$ \item $\beta
g_{x_1} + g_u \geq 0,$ on $B$, for all $|\beta| \leq 1.$
\end{enumerate}\vspace{0.01in}

\noindent We localize on the box $C_{(\alpha,M/2)}.$ Denote by
$\Omega$ be the intersection of $C_{(\alpha,M/2)}^+(u)$ with the
set $S:= \{(x_1,x_2) : (x_1,u(x_1,x_2)) \in T\}$. Notice that, in
view of $(2)$ in Lemma \ref{bound}, by choosing the width $a$ and
height $h$ of the trapezoid $T$ sufficiently small, we can
guarantee that $\Omega \subset C_{(\alpha,M/2)}.$

Define

$$H(x) = G(x)e^{x_2}\log{\left(\frac{|u_1|}{u_2}\right)},$$with

$$G(x)=g(x_1,u(x)).$$ Let

$$H(\overline{x}) = \max_{\overline{\Omega}} H(x),$$\vspace{0.01in}

\noindent and assume by contradiction that $H(\overline{x}) \geq N
\geq 1,$ for some large constant $N$ to be chosen later. Hence,

\begin{equation}\label{lowerb1}\frac{|u_1|}{u_2} \geq
\log{\left(\frac{|u_1|}{u_2}\right)} = \frac{H}{G e^{x_2}} \geq C
N, \ \ \text{at $\overline{x},$}\end{equation}\vspace{0.01in}

\noindent and also, using (1) in Lemma \ref{bound},

\begin{equation}\label{u2G}\frac{u_2}{G} \leq  C\frac{|u_1|}{H} \leq
 \frac{C'}{N}, \ \ \text{at $\overline{x}.$}\end{equation}\vspace{0.01in}

\noindent
 Furthermore,
either of the following two possibilities holds:
\begin{enumerate}
\item $\overline{x} \in F(u), $
\item $\overline{x} \in \Omega.$\end{enumerate} Let us start by showing that $(1)$ cannot occur.
 Indeed, in this case, we would have

 $$(\partial_{\nu}
\log{|H|})(\overline{x}) \leq 0,
$$ \vspace{0.01in}

\noindent where $\nu=(u_1(\overline{x}),u_2(\overline{x}))$
denotes the inner normal direction to $F(u)$ at $\overline{x}$.
Hence, at $\overline{x}$, we would have

\begin{eqnarray}\label{negnormal}
\nonumber u_1\left\{\frac{G_1}{G} +  \frac{1}{ \log{
\left(\frac{|u_1|}{u_2}\right) } }\left[\frac{u_{11}}{u_1} -
\frac{u_{21}}{u_2}\right]\right\} + u_2\left\{\frac{G_2}{G} + 1 +
\frac{1}{ \log{
\left(\frac{|u_1|}{u_2}\right) } }\left[\frac{u_{12}}{u_1} - \frac{u_{22}}{u_2}\right]\right\}=\\
\nonumber \ \\
u_1\frac{G_1}{G} + u_2\frac{G_2}{G} + u_2 + \frac{1}{ \log{
\left(\frac{|u_1|}{u_2}\right) } } \left[u_{11} -u_{22} +
u_{12}\left(\frac{u_2}{u_1} - \frac{u_1}{u_2}\right)\right] \leq
0.
\end{eqnarray}\vspace{0.01in}

\noindent Let us show that the quantity in the square bracket is
zero along $F(u)$ . The free boundary condition says that

\begin{equation}\label{fbc}u_1^2+u_2^2 = 1 \ \ \ \text{on} \ \ F(u).
\end{equation} \vspace{0.01in}

\noindent Thus, differentiating this condition along the
tangential direction $(u_2,-u_1)$ we obtain,
\begin{equation}
u_2(u_1u_{11} + u_2u_{21}) - u_1(u_1u_{12}+u_2u_{22}) = 0 \ \
\text{on $F(u)$.}
\end{equation} \vspace{0.01in}

\noindent Hence,  we deduce that

\begin{equation}
u_1u_2(u_{11} - u_{22}) + u_{12}(u_2^2-u_1^2) = 0 \ \ \text{on
$F(u)$.}
\end{equation}\vspace{0.01in}

\noindent Therefore,

$$ \left[u_{11} -u_{22} +
u_{12}\left(\frac{u_2}{u_1} -
\frac{u_1}{u_2}\right)\right](\overline{x}) = 0,$$\vspace{0.01in}

\noindent and \eqref{negnormal} reads,

\begin{equation}\label{negnormal2n}
(u_1\frac{G_1}{G} + u_2\frac{G_2}{G} + u_2)(\overline{x}) \leq 0.
\end{equation} \vspace{0.01in}

\noindent On the other hand,

\begin{equation}\label{negnormal3n}
(u_1\frac{G_1}{G} + u_2\frac{G_2}{G})(\overline{x}) =
(\frac{u_1[g_{x_1} + g_u u_1] + g_uu_2^2 }{G})(\overline{x}) =
\frac{(u_1g_{x_1} + g_u)(\overline{x})}{G(\overline{x})} \geq 0,
\end{equation}\vspace{0.01in}

\noindent according to property (iii) in the definition of $g$,and
the free boundary condition \eqref{fbc}. The inequality
\eqref{negnormal3n} together with the fact that $u_2 >0$,
contradicts \eqref{negnormal2n}.

\

\noindent \textbf{Remark.} We remark that this argument is
independent of the particular equation which is satisfied by $u$
in its positive phase. Moreover, it is easily generalized to
higher dimensions.

\

Now, we proceed to showing that by choosing $N$ sufficiently
large, we obtain a contradiction also in case $(2)$. In this case
we would have,

\begin{equation}\label{zerograd1}(\partial_i\log{|H|})(\overline{x})= 0, \ \
i=1,2,\end{equation}and

\begin{equation}\label{maxnegl}\Delta(\log|H|)(\overline{x}) \leq
0.\end{equation} For brevity, we denote by

$$L = \left(\frac{1}{ \log{ \left(\frac{|u_1|}{u_2}\right) }
}\right)(\overline{x}),$$ \vspace{0.01in}

\noindent hence, according to \eqref{lowerb1}

\begin{equation}\label{LN1} L \leq C/N.\end{equation} Then, \eqref{zerograd1} reads,

\begin{equation}\label{deriv1} \left(\frac{G_i}{G}+ \delta_{i2}
+L\left[\frac{u_{1i}}{u_1} - \frac{u_{2i}}{u_2}\right]
\right)(\overline{x}) =0 , \ \ i=1,2.
\end{equation}\vspace{0.01in}

\noindent In order to use \eqref{maxnegl}, let us compute,

\begin{eqnarray}\label{deriv2}
\partial_{ii} (\log|H|)= \frac{G_{ii}}{G} -
\frac{G_i^2}{G^2} + L\left[\frac{u_{1ii}}{u_1} -
\frac{u_{1i}^2}{u_1^2} - \frac{u_{2ii}}{u_2} +
\frac{u_{2i}^2}{u_2^2}\right]\\\nonumber \ \\\nonumber -L^2
\left[\frac{u_{1i}}{u_1} - \frac{u_{2i}}{u_2}\right]^2, \ \ \
i=1,2.
\end{eqnarray}\vspace{0.01in}

\noindent Thus, according to \eqref{maxnegl}, using that $u_1,u_2$
are harmonic functions, we get

\begin{eqnarray}\label{derivl}
\frac{\Delta G}{G} - \frac{|\nabla G|^2}{G^2} + L\left[ -
\frac{|\nabla u_1|^2}{u_1^2}+ \frac{|\nabla u_2|^2
}{u_2^2}\right]\\\nonumber \ \\\nonumber -L^2
\sum_{i=1}^2\left[\frac{u_{1i}}{u_1} - \frac{u_{2i}}{u_2}\right]^2
\leq 0, \ \ \ \text{at} \ \ \ \overline{x}.
\end{eqnarray}\vspace{0.01in}

\noindent We wish to prove that if $N$ is large enough, the
quantity $L\frac{|\nabla u_2|^2 }{u_2^2}$ is very large and hence
it dominates all the  summands in \eqref{derivl}. Toward this aim,
let us start by proving that if $N$ is sufficiently large, then

\begin{equation}\label{G2G1}
\left|\frac{G_2}{G}(\overline{x})\right| \leq \frac{1}{2},
\end{equation} \vspace{0.01in}

\noindent which combined with \eqref{deriv1} when $i=2$, implies:

\begin{equation}\label{like11} \frac{1}{2} \leq L\left|\frac{u_{12}}{u_1}
- \frac{u_{22}}{u_2}\right|\leq \frac{3}{2}.
\end{equation}\vspace{0.01in}

\noindent Indeed, from the definition of $G$ and \eqref{u2G} we
obtain immediately,

\begin{equation}\label{lap} \left|\frac{G_2}{G}(\overline{x})\right| =
\frac{|g_u| u_2}{G}(\overline{x}) \leq \frac{C}{N}.
\end{equation} \vspace{0.01in}

\noindent Therefore, for $N$ large enough, \eqref{G2G1}, hence
\eqref{like11} hold. Since $|u_1|(\overline{x}) \geq
u_2(\overline{x})$, we deduce immediately from \eqref{like11},

\begin{equation}\label{imp}
\frac{|\nabla u_2|^2}{u_2^2} \geq \frac{u_{21}^2}{u_1^2} +
\frac{u_{22}^2}{u_2^2} \geq \frac{1}{2}\left[\frac{u_{12}}{u_1} -
\frac{u_{22}}{u_2}\right]^2 \geq \frac{C}{L^2},
\end{equation}\vspace{0.01in}

\noindent which together with \eqref{LN1} gives that
$L\frac{|\nabla u_2|^2 }{u_2^2}$ is very large, for $N$ large. In
particular, according to \eqref{lap}, for $N$ large we have,

\begin{equation}\label{lap2}
\left|\frac{G_2}{G}(\overline{x})\right|^2 \leq L^2\frac{|\nabla
u_2|^2}{u_2^2}.
\end{equation}\vspace{0.01in}

\noindent Moreover, \eqref{deriv1} for $i=1$ implies

\begin{eqnarray}\label{G1b2}
\left|\frac{G_1}{G}\right|^2 = L^2 \left[\frac{u_{11}}{u_1} -
\frac{u_{21}}{u_2}\right]^2
 \leq 2L^2\left(\frac{u_{11}^2}{u_1^2} + \frac{u_{21}^2}{u_2^2}\right) \leq 2 L^2\frac{|\nabla u_2|^2}{u^2_2},
\end{eqnarray}\vspace{0.01in}

\noindent where in order to obtain \eqref{G1b2} we have used that
$ |u_1|(\overline{x}) \geq u_2(\overline{x})$, together the fact
that since $u$ is a solution to $\Delta u=0$, then
$u_{11}^2=u_{22}^2,$ and in particular

\begin{equation}\label{equal}|\nabla u_1|^2=|\nabla
u_2|^2.\end{equation} \vspace{0.01in}

\noindent Thus, combining \eqref{lap2} and \eqref{G1b2} we get,

\begin{eqnarray}\label{gradG}
\frac{|\nabla G|^2}{G^2} \leq 3 L^2\frac{|\nabla u_2|^2}{u^2_2},
\end{eqnarray}\vspace{0.01in}

\noindent and combining \eqref{imp} and \eqref{G1b2} we also get,

\begin{equation}\label{lapl}L^2\left(
\left[\frac{u_{11}}{u_1} - \frac{u_{21}}{u_2}\right]^2+
\left[\frac{u_{12}}{u_1} - \frac{u_{22}}{u_2}\right]^2\right) \leq
3 L^2\frac{|\nabla u_2|^2}{u^2_2}.
\end{equation}\vspace{0.01in}

\noindent Now, combining \eqref{derivl} with
\eqref{lowerb1},\eqref{equal},\eqref{gradG}, and \eqref{lapl},
 we obtain

\begin{eqnarray}
\frac{\Delta G}{G} -  6L^2\frac{|\nabla u_2|^2}{u^2_2}+
L\frac{|\nabla u_2|^2}{u^2_2}\left[- \frac{1}{N^2}+1\right]\leq 0,
\ \ \text{at} \ \ \overline{x}.
\end{eqnarray} \vspace{0.01in}

\noindent Therefore, for $L$ sufficiently small, that is $N$
sufficiently large, we get

\begin{eqnarray}\label{final2l}
\frac{\Delta G}{G}  + \frac{L}{2}\frac{|\nabla u_2|^2}{u^2_2} \leq
0, \ \ \text{at} \ \ \overline{x}.
\end{eqnarray}\vspace{0.01in}

\noindent On the other hand, since $\Delta u=0$, we have

\begin{eqnarray}\label{DeltaG} \frac{\Delta G}{G} =
\frac{\left[g_{x_1x_1} + 2g_{x_1u}u_1 + g_{uu}|\nabla u|^2
\right]}{G}\geq -\frac{C}{G}.
\end{eqnarray}\vspace{0.01in}

\noindent Moreover, since $H(\overline{x}) \geq N$, we have

\begin{eqnarray}\label{DeltaGf} \frac{\Delta G}{G}(\overline{x})
\geq -\frac{C}{G(\overline{x})} \geq -\frac{C}{LN}.
\end{eqnarray}\vspace{0.01in}

\noindent Therefore, combining \eqref{imp}, \eqref{final2l} and
\eqref{DeltaGf} we obtain
\begin{eqnarray}\label{final4}
-\frac{C}{LN} + \frac{C'}{L}\leq 0,
\end{eqnarray}\vspace{0.01in}

\noindent and we reach a contradiction for $N$ large.\qed

\section{Proof of Theorem \ref{main}. Non-linear operators.}

The proof follows the lines of the case when $\mathcal{F}(D^2 u)=
\Delta u. $ Precisely, with the same notation as in Section 3, we
assume by contradiction that $H(\overline{x}) \geq N \geq 1,$ for
some large constant $N$ to be chosen later. Hence, the following
three bounds hold

\begin{equation}\label{lowerbn}\frac{|u_1|}{u_2}(\overline{x}) \geq
C N,\end{equation}

\begin{equation}\label{u2Gn }\frac{u_2}{G}(\overline{x}) \leq
\frac{C}{N},\end{equation}

\begin{equation}\label{1Gn}\frac{1}{G}(\overline{x}) \leq
\frac{C}{LN}.\end{equation}\vspace{0.01in}

\noindent
 Furthermore, according to the
argument in Section 2 (and the remark following it), the maximum
must be achieved in the interior, that is $\overline{x} \in
\Omega.$

Now, we proceed to showing that by choosing $N$ sufficiently
large, we obtain a contradiction. Since $H$ achieves a maximum at
$\overline{x}$ we have,

\begin{equation}\label{zerograd}(\partial_i\log{|H|})(\overline{x})= 0, \ \
i=1,2,\end{equation}and

\begin{equation}\label{maxneg}\mathcal{L}(\log|H|)(\overline{x}) \leq
0,\end{equation} where

$$\mathcal{L}(v) = \sum_{i,j=1}^{2}a_{ij}v_{ij},$$\vspace{0.01in}

\noindent is the linearized operator associated to
$\mathcal{F}(D^2v)$. Again, we denote by

$$L = \left(\frac{1}{ \log{ \left(\frac{|u_1|}{u_2}\right) }
}\right)(\overline{x}).$$ \vspace{0.01in}

\noindent In order to use \eqref{maxneg}, let us compute,

\begin{eqnarray}\label{deriv2}
\partial_{ij} (\log|H|)= \frac{G_{ij}}{G} -
\frac{G_i G_j}{G^2} + L\left[\frac{u_{1ij}}{u_1} -
\frac{u_{1j}u_{1i}}{u_1^2} - \frac{u_{2ij}}{u_2} +
\frac{u_{2j}u_{2i}}{u_2^2}\right]\\\nonumber \ \\\nonumber -L^2
\partial_j\left(\log{\frac{|u_1|}{u_2}}\right)\partial_i\left(\log{\frac{|u_1|}{u_2}}\right),
\ \ \ i,j=1,2.
\end{eqnarray}\vspace{0.01in}

\noindent Thus, \eqref{maxneg} reads,
\begin{eqnarray}
\frac{\mathcal{L}G}{G} - \frac{1}{G^2}\sum_{i,j=1}^{2}a_{ij}G_iG_j
+ L\left[\sum_{i,j=1}^{2}a_{ij}\left(\frac{u_{2j}u_{2i}}{u_2^2} -
\frac{u_{1j}u_{1i}}{u_1^2}\right)\right] \\\nonumber \ \\\nonumber
-L^2\sum_{i,j=1}^{2}a_{ij}\partial_j\left(\log{\frac{|u_1|}{u_2}}\right)
\partial_i\left(\log{\frac{|u_1|}{u_2}}\right) \leq 0,
\ \ \text{at} \ \ \overline{x},
\end{eqnarray} \vspace{0.01in}

\noindent where we have used that $u_1,u_2$ are solutions to the
linearized equation $\mathcal{L}v=0.$ Then, by the uniform
ellipticity of $\mathcal{L}$ we derive the following inequality

\begin{eqnarray}\label{final}
\frac{\mathcal{L}G}{G} - \frac{\Lambda}{G^2}|\nabla G|^2 +
L\left[\lambda\frac{|\nabla u_2|^2}{u^2_2}-\Lambda\frac{|\nabla
u_1 |^2}{u_1^2}\right]\nonumber \ \\
\\\nonumber -L^2\Lambda\left[\left(\frac{u_{11}}{u_1}
-\frac{u_{21}}{u_2}\right)^2 + \left(\frac{u_{12}}{u_1} -
\frac{u_{22}}{u_2}\right)^2\right] \leq 0, \ \ \text{at} \ \
\overline{x}.
\end{eqnarray}\vspace{0.01in}

\noindent Again, we wish to prove that the quantity
$L\frac{|\nabla u_2|^2}{u^2_2}$ is very large, and it dominates
all the negative summands in \eqref{final}. The same argument as
in Section 3 gives that

\begin{equation}\label{imp2}
\frac{|\nabla u_2|^2}{u_2^2} \geq  \frac{C}{L^2}.
\end{equation}\vspace{0.01in}

\noindent Moreover, although \eqref{equal} is no longer valid, we
know that $u$ is a solution to $\mathcal{F}(D^2u)=0$, and
$\mathcal{F}(0)=0$. Therefore $u$ solves a linear equation with
uniformly bounded coefficients, and we get

\begin{equation}\label{compu11u22}|\nabla u_1|^2 \leq C |\nabla
u_2|^2,\end{equation}\vspace{0.01in}

\noindent for some constant $C$ depending on the ellipticity
constants $\lambda, \Lambda.$  Thus, we conclude as in the
previous section, that the following two bounds hold:

\begin{eqnarray}\label{gradGn}
\frac{|\nabla G|^2}{G^2} \leq C L^2\frac{|\nabla u_2|^2}{u^2_2},
\end{eqnarray}

\begin{equation}\label{lap3n}L^2\left(
\left[\frac{u_{11}}{u_1} - \frac{u_{21}}{u_2}\right]^2+
\left[\frac{u_{12}}{u_1} - \frac{u_{22}}{u_2}\right]^2\right) \leq
C L^2\frac{|\nabla u_2|^2}{u^2_2}.
\end{equation}\vspace{0.01in}

\noindent Combining \eqref{final} with \eqref{lowerbn},
\eqref{gradGn}, \eqref{lap3n}, we obtain

\begin{eqnarray}
\frac{\mathcal{L}G}{G} - C L^2\frac{|\nabla u_2|^2}{u^2_2} +
L\left[\lambda\frac{|\nabla u_2|^2}{u^2_2}-C \Lambda\frac{|\nabla
u_2 |^2}{N^2 u_2^2}\right] \leq 0, \ \ \text{at} \ \ \overline{x}.
\end{eqnarray}\vspace{0.01in}

\noindent Hence, for $N$ sufficiently large, that is $L$ small
enough,

\begin{eqnarray}\label{final2}
\frac{\mathcal{L}G}{G}  + C L\frac{|\nabla u_2|^2}{u^2_2} \leq 0,
\ \ \text{at} \ \ \overline{x}.
\end{eqnarray}\vspace{0.01in}

\noindent We now proceed to estimate $\mathcal{L}G/G$. Using
\eqref{u2Gn }, \eqref{1Gn} and \eqref{compu11u22}, we get
\begin{eqnarray}\label{LGn} \frac{|\mathcal{L}G|}{G} =
\frac{1}{G}|\sum_{i,j=1}^{2}a_{ij}\left[g_{x_1x_1}\delta_{i1}\delta_{j1}
+ 2g_{x_1u}\delta_{i1}u_j + g_{uu}u_iu_j +
g_uu_{ij}\right]|\\\nonumber \ \\\nonumber \leq \frac{1}{G}(C + C
|\nabla u_2|) \leq \frac{C}{LN} + \epsilon L \frac{|\nabla
u_2|^2}{u^2_2} + \frac{1}{4\epsilon L }\frac{u_2^2}{G^2}
\\\nonumber \ \\\nonumber \leq \frac{C}{LN} + \epsilon L \frac{|\nabla
u_2|^2}{u^2_2} + \frac{1}{4\epsilon L }\frac{C}{N^2}.
\end{eqnarray}\vspace{0.01in}

\noindent Hence, for $\epsilon$ small, combining \eqref{final2}
with \eqref{imp2}, and \eqref{LGn} we get,

\begin{eqnarray}\label{final11} -\frac{C}{LN} -\frac{C}{LN^2}
+\frac{C}{L} \leq 0,\end{eqnarray}\vspace{0.01in}

\noindent that is a contradiction for $N$ large enough.

\section{The $p$-Laplace operator.}

In this section, we generalize the a-priori bound in Theorem
\ref{main}, to the case when  $u$ is a classical solution to the
following one-phase free boundary problem:

\begin{equation}\label{fbproblemp} \left \{
\begin{array}{ll}
     \text{div}(|\nabla u|^{p-2}\nabla u)=0,   & \hbox{in $C^+_{M}(u)$,} \\
    |\nabla u| =1, & \hbox{on $F(u)$,} \\
   u_2 >0 , & \hbox{on $\overline{C^+_{M}(u)}$,} \\
\end{array}\right.
\end{equation}\vspace{0.01in}

\noindent with $1<p < \infty.$ The setting will be the same as in
Section 1, that is we assume that $0 \in F(u)$, and

$$B_{2\alpha }(0, M/2) \Subset
 C^+_{M}(u), \  \ \text{and} \ \
B_{2\alpha }(0, -M/2) \Subset C^-_{M}(u),$$\vspace{0.01in}

\noindent
 for some constant
$\alpha <1/2.$

Here and henceforth, $C, C'$ will denote constants depending
possibly on $M, p$ and $\alpha.$ The following technical Lemma
still holds.

\begin{lem}\label{boundp} Let $u$ be a solution to \eqref{fbproblemp}. Then, there exist constants $C,C'$, such that
\begin{enumerate}\item $|\nabla u| \leq C'$ on
$C_{(\alpha,M/2)},$ \item $u \geq C$ on
$\mathcal{B}_{\alpha}(0,M/2).$ \end{enumerate}
\end{lem}

We wish to prove the following result.

\begin{thm}\label{mainp} There exists constants C, s,
depending on $M, \lambda, \Lambda, \alpha, p$ such that

$$\frac{|u_1|}{u_2}\leq C \ \ \textit{on} \ \ B_s(0).$$ \vspace{0.01in}

\noindent In particular $Lip(s) \leq C.$
\end{thm}
\begin{proof} The proof follows the lines of the proof of Theorem
\ref{main}. Again, we introduce the function $H$ and we assume
that it achieves a large maximum at $\overline{x}$, which
according to the argument in Section 2 must be an interior point,
i.e. $\overline{x} \in \Omega.$ Then,

\begin{equation}\label{zerogradp}(\partial_i\log{|H|})(\overline{x})= 0, \ \
i=1,2,\end{equation}and

\begin{equation}\label{maxnegp}\mathcal{L}(\log|H|)(\overline{x}) \leq
0,\end{equation} where

$$\mathcal{L}(v) = \sum_{i,j=1}^{2}a_{ij}v_{ij}, \ \ a_{ij} = \delta_{ij} + (p-2)\frac{v_iv_j}{|\nabla
v|^2}.$$\vspace{0.01in}

\noindent In particular, $\mathcal{L}$ is uniformly elliptic with
constants $\lambda = \min\{1,p-1\}, \Lambda = 1 + |p-2|.$ Notice
that since $|\nabla u|>0,$ $u$ is a solution to $\mathcal{L}v=0.$
Again, we denote by

$$L = \left(\frac{1}{ \log{ \left(\frac{|u_1|}{u_2}\right) }
}\right)(\overline{x}).$$\vspace{0.01in}

\noindent According to \eqref{maxnegp}, from formula
\eqref{deriv2} and from the uniform ellipticity of $\mathcal{L}$,
we derive the following inequality

\begin{eqnarray*}\label{finalp}
\frac{\mathcal{L}G}{G} - \frac{\Lambda}{G^2}|\nabla G|^2 +
L\left[\lambda\frac{|\nabla u_2|^2}{u^2_2}-\Lambda\frac{|\nabla
u_1 |^2}{u_1^2}\right] +
L\left[\sum_{i,j=1}^{2}a_{ij}\left[\frac{u_{1ij}}{u_1}  -
\frac{u_{2ij}}{u_2}\right]\right]
\\\nonumber \ \\\nonumber -L^2\Lambda\left[\left(\frac{u_{11}}{u_1}
-\frac{u_{21}}{u_2}\right)^2 + \left(\frac{u_{12}}{u_1} -
\frac{u_{22}}{u_2}\right)^2\right] \leq 0, \ \ \text{at} \ \
\overline{x}.
\end{eqnarray*}\vspace{0.01in}

\noindent The difference between this inequality and
\eqref{final}, consists in the presence on the term
$L\left[\sum_{i,j=1}^{2}a_{ij}\left[\frac{u_{1ij}}{u_1}  -
\frac{u_{2ij}}{u_2}\right]\right]$, which appears since $u_1$ and
$u_2$ are not solutions to $\mathcal{L}v=0.$ Thus, if we show that

\begin{equation}\label{diffp}
\left|\sum_{i,j=1}^{2}a_{ij}\left[\frac{u_{1ij}}{u_1}  -
\frac{u_{2ij}}{u_2}\right]\right| \leq
\frac{\lambda}{2}\frac{|\nabla u_2|^2}{u^2_2},
\end{equation}\vspace{0.01in}

\noindent we obtain that
\begin{eqnarray}\label{finalpp}
\frac{\mathcal{L}G}{G} - \frac{\Lambda}{G^2}|\nabla G|^2 +
L\left[\frac{\lambda}{2}\frac{|\nabla
u_2|^2}{u^2_2}-\Lambda\frac{|\nabla u_1 |^2}{u_1^2}\right]
\\\nonumber \ \\\nonumber -L^2\Lambda\left[\left(\frac{u_{11}}{u_1}
-\frac{u_{21}}{u_2}\right)^2 + \left(\frac{u_{12}}{u_1} -
\frac{u_{22}}{u_2}\right)^2\right] \leq 0, \ \ \text{at} \ \
\overline{x},
\end{eqnarray} \vspace{0.01in}

\noindent and we reach a contradiction as in the previous section.

In order to prove \eqref{diffp}, we start by differentiating the
equation $\mathcal{L}u=0.$ We get,

\begin{equation}
0=\partial_k(\mathcal{L}u) = \sum_{i,j=1}^2 a_{ij}u_{ijk} +
\sum_{i,j=1}^2 \partial_k (a_{ij})u_{ij},
\end{equation}\vspace{0.01in}

\noindent hence

\begin{equation}
\label{deriv3} \sum_{i,j=1}^2 a_{ij}u_{ijk} + 2(p-2)\sum_{i,j=1}^2
\left(\frac{u_iu_{jk}}{|\nabla u |^2} - u_iu_j\frac{\sum_{l=1}^2
u_lu_{lk}}{|\nabla u|^4} \right)u_{ij}=0.
\end{equation}\vspace{0.01in}

\noindent Now, using \eqref{deriv3}, set $\epsilon =
\frac{\lambda}{16(p-1)}$ we obtain that

\begin{eqnarray*} \left|\frac{1}{u_1}\sum_{i,j=1}^2 a_{ij}u_{1ij}
- \frac{1}{u_2}\sum_{i,j=1}^2 a_{ij}u_{2ij} \right|= \\\nonumber \
\\\nonumber 2\frac{p-1}{|\nabla u|^2}\left|u_{11}^2 - u_{22}^2 +
u_{12}\left(\frac{u_2}{u_1}-\frac{u_1}{u_2}\right)(u_{11}+u_{22})
\right| \leq \\\nonumber \ \\ 2\frac{p-1}{|\nabla
u|^2}\left[u_{11}^2 + u_{22}^2 + 2\epsilon
u^2_{12}\left(\frac{u_2^2}{u_1^2}+\frac{u_1^2}{u_2^2}\right) +
\frac{1}{2\epsilon}(u_{11}^2+u_{22}^2) \right]\leq
\\\nonumber \ \\\label{beflast}
C\frac{1}{\epsilon}(p-1)\frac{|\nabla u_2|^2}{N^2u_2^2}+4\epsilon
\frac{p-1}{u_2^2} |\nabla u_2|^2\leq\\\nonumber \ \\
\frac{\lambda}{2} \frac{|\nabla u_2|^2}{u_2^2}
\end{eqnarray*}\vspace{0.01in}

\noindent as long as $N$ is large enough.
\end{proof}

\section{A priori bound for 3D free boundary graphs.}

In this section we extend our result in 3D, for the case when $u$
is harmonic in its positive phase. We intend to highlight the
difficulties which arise when trying to adapt our technique to
higher dimensions.

Assume that $u$ is a classical solution to the following one-phase
free boundary problem:

\begin{equation}\label{fbproblem2} \left \{
\begin{array}{ll}
     \Delta u=0,   & \hbox{in $C^+_{M}(u)$,} \\
    |\nabla u| =1, & \hbox{on $F(u)$,} \\
   u_3 >0 , & \hbox{on $\overline{C^+_{M}(u)}$.} \\
\end{array}\right.
\end{equation}\vspace{0.01in}

\noindent Furthermore, assume that $0 \in F(u)$, and that

$$B_{2\alpha }(0, M/2) \Subset
 C^+_{M}(u),  \ \ \text{and} \ \
B_{2\alpha }(0, -M/2) \Subset C^-_{M}(u),$$\vspace{0.01in}

\noindent
 for some constant
$\alpha <1/2.$

\begin{thm}\label{main3} There exists constants C, s,
depending on $M, \alpha$ such that $$Lip(s) \leq C.$$
\end{thm}

\begin{proof}Let $g=g(r,u)$ be the function introduced in Section
3. One can easily construct such function, so that it satisfies
the following condition:\vspace{0.01in}

\noindent

(iv) $\ \ |\nabla g|^2 \leq C g  \ \ \text{near} \ \ \partial T
\setminus B.$

\vspace{0.1in}

\noindent We localize on the box $C_{(\alpha,M/2)}.$ Denote by
$\Omega$ be the intersection of $C_{(\alpha,M/2)}^+(u)$ with the
set $S:= \{x : (|x'|,u(x)) \in T\}$. Again, in view of $(2)$ in
Lemma \ref{bound},  $\Omega \subset C_{(\alpha,M/2)}.$

Define

$$H(x) =
G(x)e^{x_3}\log{\left(\frac{|\nabla_{x'}u|}{u_3}\right)},$$ with

$$G(x)=g(|x'|,u(x)).$$ Let

$$H(\overline{x}) = \max_{\overline{\Omega}} H(x),$$\vspace{0.01in}

\noindent and assume by contradiction that $H(\overline{x}) \geq N
\geq 1,$ for some large constant $N$ to be chosen later. Without
loss of generality we can assume that
$|\nabla_{x'}u|(\overline{x})=u_1(\overline{x}),$ hence in
particular

\begin{equation}\partial_i(|\nabla_{x'}u|)(\overline{x})=u_{1i}(\overline{x}).\end{equation}Also,

\begin{equation}\label{lowerb13}\frac{u_1}{u_3} \geq
\frac{H}{G e^{x_3}} \geq C N, \ \ \text{at
$\overline{x}.$}\end{equation}\vspace{0.01in}

\noindent By the same argument as in 2D, one can deduce that
$\overline{x}$ is an interior point.

We proceed to showing that by choosing $N$ sufficiently large, we
obtain a contradiction. We have,

\begin{equation}\label{zerograd13}(\partial_i\log{|H|})(\overline{x})= 0, \ \
i=1,2,3\end{equation} and

\begin{equation}\label{maxnegl3}\Delta(\log|H|)(\overline{x}) \leq
0.\end{equation}For brevity, we denote by

$$L = \left(\frac{1}{ \log{ \left(\frac{u_1}{u_3}\right) }
}\right)(\overline{x}),$$\vspace{0.01in}

\noindent hence, according to \eqref{lowerb13}

\begin{equation}\label{LN} L \leq C/N.\end{equation} Then, \eqref{zerograd13} reads,

\begin{equation}\label{deriv13} \left(\frac{G_i}{G}+ \delta_{i3}
+L\left[\frac{u_{1i}}{u_1} - \frac{u_{3i}}{u_3}\right]
\right)(\overline{x}) =0 , \ \ i=1,2,3.
\end{equation} \vspace{0.01in}

\noindent In particular, at $\overline{x},$

\begin{equation}\label{new3}
L^2 \sum_{i=1}^3\left[\frac{u_{1i}}{u_1} -
\frac{u_{3i}}{u_3}\right]^2 \leq \frac{|\nabla G|^2}{G^2} + 1 .
\end{equation}\vspace{0.01in}

\noindent In order to use \eqref{maxnegl3}, let us compute at
$\overline{x}$,

\begin{eqnarray*}\label{deriv23}
\partial_{ii} (\log|H|)= \frac{G_{ii}}{G} -
\frac{G_i^2}{G^2} +
L\left[\frac{\partial_{ii}|\nabla_{x'}u|}{|\nabla_{x'}u|} -
\frac{u_{1i}^2}{u_1^2} - \frac{u_{3ii}}{u_3} +
\frac{u_{3i}^2}{u_3^2}\right]\\\nonumber \ \\\nonumber -L^2
\left[\frac{u_{1i}}{u_1} - \frac{u_{3i}}{u_3}\right]^2, \ \
i=1,2,3,
\end{eqnarray*}\vspace{0.01in}

\noindent and,

$$\frac{\Delta|\nabla_{x'}u|}{|\nabla_{x'}u|}= \sum_{i=1}^{3}\frac{u_{2i}^2}{u_1^2}.  $$ \vspace{0.01in}

\noindent Thus, according to \eqref{maxnegl3}, using \eqref{new3}
together with the fact that $u_3$ is harmonic, we get

\begin{eqnarray}\label{derivl3}
\frac{\Delta G}{G} - 2\frac{|\nabla G|^2}{G^2} -1 +
L\left[\sum_{i=1}^{3}\frac{(u_{2i}^2-u_{1i}^2)}{u_1^2} +
\frac{|\nabla u_3|^2 }{u_3^2}\right] \leq 0, \ \ \ \text{at} \ \ \
\overline{x}.
\end{eqnarray}\vspace{0.01in}

\noindent With similar computations as in 2D, one has that
$L\frac{|\nabla u_3|^2 }{u_3^2}$ is very large, and it dominates
all the other summands. Precisely

\begin{equation}\label{imp3d}
\frac{|\nabla u_3|^2 }{u_3^2} \geq
\frac{1}{2}\left[\frac{u_{13}}{u_1} - \frac{u_{33}}{u_3}\right]^2
\geq \frac{C}{L^2}.
\end{equation}\vspace{0.01in}

\noindent Now, from \eqref{derivl3}, using that $u$ is harmonic,
we obtain at $\overline{x},$

\begin{eqnarray}\label{derivl3n} \ \ \ \ \ \ \ \ \
\frac{\Delta G}{G} - 2\frac{|\nabla G|^2}{G^2} -1 + L\left[
-\frac{u_{13}^2}{u_1^2} + \frac{u_{33}^2}{u_1^2} +
2\frac{u_{11}u_{33}}{u^2_1} + \frac{u_{31}^2}{u_3^2} +
\frac{u_{33}^2}{u_3^2}\right] \leq 0.
\end{eqnarray}\vspace{0.01in}

\noindent On the other hand, \eqref{deriv13} for $i=1$ gives,

\begin{equation}L\left|\frac{u_{11}}{u_1}\right| \leq
L\left|\frac{u_{31}}{u_3}\right|+  \left|\frac{G_1}{G}\right|\ \ \
\text{at} \ \ \ \overline{x}.
\end{equation}\vspace{0.01in}

\noindent Hence, using that $u_1/u_3 \geq 1/L$ at $\overline{x}$,
we have

\begin{equation} 2L\left|\frac{u_{11}u_{33}}{u^2_1}\right| \leq L
\left[\epsilon \frac{u_{31}^2}{u_3^2} +
\frac{1}{\epsilon}\frac{u_{33}^2}{u_1^2}\right] +
\left[\frac{1}{\epsilon}\frac{G_1^2}{G^2} + \epsilon
L\frac{u_{33}^2}{u_3^2}\right]\ \ \ \text{at} \ \ \ \overline{x}.
\end{equation}\vspace{0.01in}

\noindent Combining this estimate with \eqref{derivl3n} we obtain,
$(\epsilon=1/2$)

\begin{equation}\label{3d}
\frac{\Delta G}{G} - 4\frac{|\nabla G|^2}{G^2} -1 + L\left[
-\frac{u_{13}^2}{u_1^2} - \frac{u_{33}^2}{u_1^2} + \frac{1}{2}
\frac{u_{31}^2}{u_3^2} + \frac{1}{2}\frac{u_{33}^2}{u_3^2}\right]
\leq 0, \ \  \text{at} \ \  \overline{x}.
\end{equation}\vspace{0.01in}

\noindent Hence, for $N$ large enough, using \eqref{imp3d} we get

\begin{eqnarray}\label{3df} \frac{\Delta G}{G} - 4\frac{|\nabla
G|^2}{G^2} -1 + \frac{C}{L} \leq 0, \ \  \text{at} \ \
\overline{x}.
\end{eqnarray}\vspace{0.01in}

\noindent We can now reach a contradiction as in the 2 dimensional
case, using that, according to property (iv) of $g$, we have
$|\nabla G|^2/G^2 \leq C/G$ at $\overline{x}.$
\end{proof}

\end{document}